\documentclass{elsart}

\usepackage{natbib}

\usepackage{amssymb}

\usepackage[francais]{babel}
\usepackage[latin1]{inputenc}
\usepackage{fancyhdr}
\usepackage{amsmath}
\usepackage{amsfonts}
\usepackage{amsthm}
\usepackage{latexsym}
\usepackage{layout}
\usepackage{hyperref}
\usepackage{layout}
\usepackage{graphicx,subfigure}

\newcommand{\1}{{\rm 1}\kern-0.26em{\rm I}}
\newcommand{\E}{\mathbb{E}}
\newcommand{\p}{{\rm I}\kern-0.18em{\rm P}}

\newcommand{\R}{{\rm I}\kern-0.20em{\rm R}}
\newcommand{\bx}{\boldsymbol{x}}

\newtheorem{theorem}{Theorem}[section]

\newtheorem{lemma}{Lemma}[section]

\newtheorem{corollary}{Corollary}[section]

\newtheorem{remark}{Remark}

\renewcommand{\theequation}{\thesection.\arabic{equation}}

\begin{document}

\begin{frontmatter}



\title{Asymptotic normality for estimators of the additive regression components under random censorship}


\author[label1]{M. Debbarh}
\ead{debbarh@ccr.jussieu.fr}
\author[label1],label2{V. Viallon}
\ead{viallon@ccr.jussieu.fr}

\address[label1]{L.S.T.A, Université Paris 6, 175, rue du Chevaleret,
75013 Paris, France}
\address[label2]{Unité de Biostatistique, Hôpital Cochin, Université Paris-Descartes, Paris, F-75014 FRANCE}

\begin{abstract}
We establish asymptotic normality for estimators of the additive
regression components under random censorship. To build our
estimators, we couple the marginal integration method \cite{Newey}
with an initial \emph{Inverse Probability of Censoring Weighted}
estimator of the multivariate censored regression function
introduced by \cite{Carbonez1995} and \cite{Kohler2002}.
Asymptotic confidence bands are derived from our result.

\end{abstract}


\end{frontmatter}

\section{Introduction}
\setcounter{equation}{0} \setcounter{theo}{0} \setcounter{lem}{0}
\setcounter{cor}{0} \setcounter{rem}{0} \setcounter{prop}{0}
\setcounter{fact}{0}

Censored data arise in many statistical application domains,
especially in epidemiology and reliability. When studying the
relationship between a depending variable and covariates,
nonparametric estimates are of particular interest in the presence
of censored data, because no scatter plots can be drawn to detect
the possibly complex form of this relationship. Several estimators
have already been proposed to estimate the regression function.
The key idea is to transform the observed data and derive
so-called \emph{synthetic data} estimators \cite{BuckleyJames}.
For instance, Fan and Gijbels developed a local version of the
parametric Buckley and James estimator, while Carbonez \emph{et
al.} and Kohler \emph{et al.} studied a nonparametric
\emph{Inverse Probability of Censoring Weighted} [I.P.C.W.]
estimator \cite{Carbonez1995,Kohler2002,BrunelComte1}. \vskip5pt

\noindent However, in most situations, the encountered variables
(typically the onset of a disease in epidemiology) are related to
numerous factors or predictors. Consider the particular case of
the construction of a disease risk prediction tool, that is, the
estimation of the probability of developing the disease, given a
number (about ten generally) of predictors. In this setting,
classical nonparametric estimators, such as kernel-type
estimators, are unsuitable because of the well-known \emph{curse
of dimensionality}: for instance, the rate of convergence for
nonparametric estimators of the conditional survival function is
increasing in the regressor dimension in a censored setting
\cite{Dabrowska}, just like it is the case in an uncensored one
(see, e.g., \cite{Hardle90} and the relevant references therein).
One common solution to get round this issue is to work under the
additive model assumption, when possible. In the uncensored case,
several methods have been proposed to estimate the additive
regression function. We shall evoke, among others, the methods
based on $B$-splines \cite{Stone85}, on the \emph{backfitting}
algorithm \cite{HastieetTibshirani} and on marginal integration
\cite{newey1994,Auestad,Linton1995}. In \cite{FanGijbels}, it is
shown that the backfitting ideas also applies to censored data.
Here, following the ideas introduced in \cite{DebbarhViallon3}, we
make use of the marginal integration method, coupled with initial
multivariate nonparametric I.P.C.W. estimators to provide an
estimator for the additive censored regression function. At this
point, it is noteworthy that the developments we propose here for
$I.P.C.W$-type estimators shall apply with minor modifications to
cope with other synthetic data estimators.

In former works, we established the uniform consistency rate
\cite{DebbarhViallon}, the mean-square convergence rate
\cite{DebbarhViallon3} and a uniform law of the logarithm
\cite{DebbarhViallon5} for such estimators of the additive
regression function in the presence of censored data. In
\cite{DebbarhViallon5}, we also proposed a method to construct
simultaneous almost certainty bands, that is confidence bands
which contain the true value of the additive component with
asymptotical probability one, uniformly over the predictor domain.
Obviously, those kinds of confidence bands may be very
conservative, and classical confidence intervals derived from some
normal law may be desirable. To construct such intervals is one of
the aims of the present paper, which is organized as follows.
After having recalled how to construct estimators for the additive
components in Section \ref{sec_notations}, we establish their
asymptotic normality in Section \ref{sec_Results}. This limit law
completes the one obtained by \cite{Sarad2000} in the uncensored
case (see also \cite{Linton1995} and \cite{Sperlich2003}). Then,
in Section \ref{sec_Applications}, we show how to obtain
confidence intervals from the aforementioned convergence in law.
Finally, Section \ref{sec_Proofs} is devoted to the proof of our
result.

\section{Notations} \label{sec_notations}
\setcounter{equation}{0} \setcounter{theo}{0} \setcounter{lem}{0}
\setcounter{cor}{0} \setcounter{rem}{0} \setcounter{prop}{0}
\setcounter{fact}{0} Consider the triple $(Y, C, \bf{X})$ defined
in $\mathbb{R}^+ \times \mathbb{R}^+ \times \mathbb{R}^d$,
$d\geq1$, where $Y$ is the variable of interest, $C$ the censoring
variable and ${\bf X}=(X_1,...,X_d)$ a vector of concomitant
variables. Throughout, we work with a sample $\{(Y_i,C_i,{\bf
X}_i)_{1\leq i \leq n}\}$ of independent and identically
distributed replicae of $(Y,C,{\bf X})$. Actually, in the right
censorship model, $Y_i$ and $C_i$ are not observed and only
$Z_i=\min\{Y_i,C_i\}$ and $\delta_i=\mathbb{I}_{\{{Y_i} \leq
C_i\}}$, $1\leq i\leq n$, are at our disposal, $\mathbb{I}_E$
standing for the indicator function of $E$. Accordingly, the
observed sample is ${\cal D}_n =\{(Z_i, \delta_i,
\boldsymbol{X}_i), i=1, \ldots, n\}$, and for all $t\in
\mathbb{R}^+$, we set $F(t)=P(Y> t)$, $G(t)=P(C> t)$ and $H(t)=
P(Z>t)$ the right continuous survival functions pertaining to $Y$,
$C$ and $Z$ respectively.\vskip5pt

\noindent Further denote by $\psi$ a given real measurable
function. In this paper, we are concerned with the regression
function of $\psi(Y)$ evaluated at $\bf{X}={\bx}$, in the
particular case where this function is additive,
\begin{eqnarray}\label{add_component}
m_{\psi}({\bx})& = & \E \left(\psi (Y) \mid \bf{X}={\bx}\right),~
\forall ~{\bx}=(x_1,...,x_d) \in \mathbb{R}^d,\mathbb{} \label{fdereg} \\
&= & \mu + \sum_{\ell=1}^d m_\ell(x_\ell).\label{additivemodel}
\end{eqnarray}
In view of (\ref{additivemodel}), the so-called \emph{additive
components} $m_\ell$, $\ell=1,...,d,$ as well as the constant
$\mu$ are only defined up to an additive constant. Therefore, it
is quite common to work under the identifiability assumption $\E
m_\ell(X_\ell)=0, \ell=1, ..., d$, which ensures that $\mu=\E
\psi(Y)$.\vskip5pt

\noindent Let $(h_n)_{n\geq 1}$ and $(h_{\ell,n})_{n\geq 1}$,
$\ell=1, \ldots, d,$ be $d+1$ sequences of real numbers and denote
by $f$ the density function of the covariate ${\bf X}$. Introduce
$\hat f_n$ the Akaike-Parzen-Rosenblatt (\cite{Akaike1954},
\cite{Parzen}, \cite{Rosemblatt}) estimator of $f$ pertaining to
$K$ and $h_n$,
\begin{eqnarray*}
\hat f_n({\bx})= \frac{1}{nh_n^d} \sum_{j=1}^n K\Big(\frac{{\bx}
-{\bf X}_j}{h_n}\Big),
\end{eqnarray*}
where $K$ is a given convolution kernel in $\mathbb{R}^d$. Further
introduce $d$ kernels $K_1,\ldots,K_d$, defined in $\mathbb{R}$.
In the sequel, we denote by $G^\star_n$ the Kaplan-Meier \cite{KM}
estimate of $G$. Namely, adopting the convention
$\prod_{\varnothing}=1$ and $0^0=1$, we have, for all $y \in
\mathbb{R}$,
\begin{equation}
G^\star_n (y) =1 -\!\!\! \prod_{1\leq i \leq
n}\!\!\!{\Big(\frac{N_n(Z_i)-1}{N_n(Z_i)}\Big)}^{\beta_i},
\label{Gnstar}
\end{equation}
with $\beta_i=\mathbb{I}_{\{Z_i\leq y\}}(1-\delta_i)$ and
$N_n(y)=\sum_{j=1}^n\mathbb{I}_{\{Z_j\geq y\}}$. \vskip5pt

\noindent To estimate the regression function defined in
(\ref{fdereg}) (at this point, we do not work under the additive
assumption (\ref{additivemodel}) yet), we propose the following
estimator (see \cite{Carbonez1995}, \cite{DebbarhViallon3},
\cite{Jones1994}, \cite{Kohler2002}  and \cite{MaillotViallon}),
\begin{equation}
\widetilde{m}_{\psi,n}^\star ({\bx}) = \sum_{i=1}^n
W_{n,i}({\bx})\frac{\delta_i \psi(Z_i)}{G^\star_n(Z_i)} ~~
\mbox{where}~~W_{n,i}({\bx})= \frac{\prod_{\ell=1}^{d}
\frac{1}{h_{\ell,n}} K_{\ell}\big({\frac{x_{\ell} -
X_{i,l}}{h_{\ell,n}}}\big)}{n\hat f_n ({\bf
X}_i)}.\label{estcarbo1}
\end{equation}
Adopting the convention $0/0=0$, $\widetilde{m}_{\psi,n}^\star$ is
properly defined since $G^\star_n(Z_i)=0$ if and only if
$Z_i=Z_{(n)}$ and $\delta_{(n)}=0$, where $Z_{(k)}$ is the $k$-th
ordered statistic associated to the sample $(Z_1,...,Z_n)$ for
$k=1,...,n$ and $\delta_{(k)}$ is the $\delta_j$ corresponding to
$Z_{(k)}=Z_j$. \vskip5pt

\noindent For all ${\bx}=(x_1,..,x_{d})\in\mathbb{R}^d$ and every
$\ell=1,...,d$, further set ${\bx}_{-\ell}=(x_1,..,$
$x_{\ell-1},x_{\ell+1},$ $..,$ $x_d)$. Under the assumption
(\ref{additivemodel}), we will estimate the additive components
via the marginal integration method (\cite{Linton1995},
\cite{Newey}). Towards this aim, we introduce $d$ given density
functions defined in $\mathbb{R}$, $q_1,...,q_d$. Then, setting
$q({\bx}) = \prod_{\ell=1}^d q_\ell(x_\ell)$ and, for
$\ell=1,...,d$, $q_{-\ell}({\bx}_{-\ell}) = \prod_{j \neq \ell}
q_j(x_j)$, consider, for $\ell=1,...,d$, the quantities
\begin{equation}
\eta_{\ell}(x_{\ell}) = \int_{\mathbb{R}^{d-1}} m_{\psi}({\bx})
q_{-\ell}({\bx}_{-\ell}) d{\bx}_{-\ell} - \int_{\mathbb{R}^d}
m_{\psi}({\bx}) q({\bx}) d{\bx}. \label{additive_component}
\end{equation}
It is straightforward that the two following equalities hold,
\begin{eqnarray}
&&\eta_\ell(x_\ell) = m_\ell(x_\ell) - \int_{\mathbb{R}} m_\ell(z)
q_\ell(z)dz, \quad \ell=1,...,d, \label{relation_additive_component}\\
&&m_{\psi}({\bx})= \sum_{\ell=1}^d \eta_\ell(x_\ell) +
\int_{\mathbb{R}^d} m_{\psi}({\bf z})q({\bf z}) d{\bf z}
\label{additive_component_marginale}.
\end{eqnarray}
In view of (\ref{relation_additive_component}) and
(\ref{additive_component_marginale}), $\eta_\ell$ and $m_\ell$ are
equal up to an additive constant, in such a way that $\eta_\ell$
is actually an additive component too, which fulfils an
alternative identifiability condition.

\begin{remark}\label{rem1}
Observe that $\eta_\ell=m_\ell$ for the particular choice
$q_\ell=f_\ell$, $f_\ell$ denoting the density function pertaining
to $X_\ell$, $\ell=1, ..., d$. However, $f_\ell$ being unknown in
most situations, we have most often $\eta_\ell\neq m_\ell$.
\end{remark}

\noindent From (\ref{estcarbo1}) and (\ref{additive_component}), a
natural estimator of the $\ell$-th additive component $\eta_\ell$
evaluated at $x_\ell$ is given, for $\ell=1,...,d$, by
\begin{equation}\label{enq1}
\widehat \eta^\star_{\ell}(x_{\ell}) = \int_{\mathbb{R}^{d-1}}\!
\widetilde{m}_{\psi,n}^{\star}({\bx}) q_{-\ell}({\bx}_{-\ell})
d{\bx}_{-\ell} - \int_{\mathbb{R}^d}\!
\widetilde{m}_{\psi,n}^{\star}({\bx}) q({\bx}) d{\bx},\
\end{equation}
from which an estimator $\widehat m_{\psi,add}^\star$ of the
additive regression function can easily be deduced (see
(\ref{additive_component_marginale})),
\begin{eqnarray}
\widehat m_{\psi,add}^\star(\bx) & = & \sum_{\ell=1}^d \hat
\eta^\star_\ell(x_\ell) + \int_{\mathbb{R}^d}
\widetilde{m}_{\psi,n}^{\star}({\bx}) q({\bx})d{\bx}
\label{estim_add}.
\end{eqnarray}

\section{Hypotheses and Results} \label{sec_Results}
\setcounter{equation}{0} \setcounter{theo}{0} \setcounter{lem}{0}
\setcounter{cor}{0} \setcounter{rem}{0} \setcounter{prop}{0}
\setcounter{fact}{0}

\noindent These preliminaries being given, we introduce the
assumptions to be made to state our result. First, consider the
hypotheses pertaining to $(Y,C,{\bf X})$.
\begin{tabbing}
$(C.1)\;\;$\=\ $C$ and $({\bf X},Y)$ are independent. \\
$(C.2)$\>\ $G$ is continuous on $\mathbb{R}^+$.\\
$(C.3)$\>\ There exists a finite constant $M$ such that
$\sup_{y}|\psi(y)| \leq M$.\\
$(C.4)$\>\ $m_\psi$ is $k$-times continuously differentiable, $ k
\geq 1$, and\\
\>\ \\[-0.5cm] \>\   \quad\quad\quad
${\sup_{\bx}\Big|\displaystyle\frac{\partial^k}{\partial
x_\ell^k}\ m_{\psi} ({\bx}) \Big|<\infty};~\ell=1,...,d$.
\end{tabbing}

\noindent As mentioned in \cite{GrossLai}, functionals of the
(conditional) law can generally not be estimated on the complete
support when the variable of interest is right-censored.
Accordingly, in order to establish our results, we will work under
the assumption $({\bf A})$ that will be said to hold if either
$({\bf A})(i)$ or $({\bf A})(ii)$ below holds. Denote by
$T_L=\sup\{t:L(t)>0\}$ the upper endpoint of the distribution of a
random variable with right continuous survival function $L$.
\begin{tabbing}
$({\bf A})(i)\;\;$ \=\ There exists
a $\tau_0<T_H$ such that  $\psi=0$ on $(\tau_0,\infty)$.\\[0.2cm]
$({\bf A})(ii)$ \>\ $(a)\;\;$ For a given $k/(2k+1)<p\leq 1/2$,
$\big|\int_0^{T_H}
 F^{-p/(1-p)}dG\big|<\infty$;\\
\>\ $(b)\;\;$ $T_F<T_G$;\\
\>\ $(c)\;\;$ $ n^{2p-1}h^{-1}_{\ell,n}|\log(h_{\ell,n})|
\rightarrow \infty$, as $n\rightarrow\infty$, for every
$\ell=1,...,d$.
\end{tabbing}

\begin{remark}
$(i)\;\;$It is noteworthy that the assumption $({\bf A})(ii)$
allows to consider the estimation of the "classical" regression
function, which corresponds to the choice $\psi(y)=y$. On the
other hand, normality for estimators of functionals such as the
conditional distribution function $\p(Y\leq \tau_0|{\bf X})$ can
be obtained under weaker conditions, when restricting ourselves to
$\tau_0<T_H$. \vskip2pt

\noindent $(ii)\;\;$When working under $({\bf A})(ii)$, assumption
$(C.3)$ can be weakened. In this setting, it is indeed sufficient
to work under $(\widetilde{C.3})$ below.
\begin{tabbing}
$(\widetilde{C.3})\;\;$ There exists a finite constant $M$ such
that $\sup_{y\leq\tau_0}|\psi(y)| \leq M$.
\end{tabbing}
\noindent $(iii)\;\;$It is also noteworthy that condition $(C.1)$
is stronger than the conditional independence of $C$ and $Y$ given
${\bf X}$, under which Beran \cite{BE81} worked to build an
estimator of the conditional survival function (see also
\cite{Dabrowska} and \cite{DD06}). Note, however, that the two
assumptions coincide if $C$ and ${\bf X}$ are independent. In
other respect, to use Beran's local Kaplan-Meier estimator, the
censoring has to be locally fair, that is $\mathbb{P}[C\geq t\mid
{\bf X}={\bx}]>0$ whenever $\mathbb{P}[Y\geq t\mid {\bf
X}={\bx}]>0$. Here, we basically only suppose that $G(t)>0$
whenever $F(t)>0$, which is, on its turn, a weaker assumption. For
a nice discussion on the difference between Beran's estimator and
Carbonez et al.'s estimator, we refer to \cite{Carbonez1995}.
\end{remark}\vskip3pt

\noindent Denote by $\mathcal{C}_1, ...,\mathcal{C}_d$, $d$
compact intervals of $\mathbb{R}$ and set
$\mathcal{C}=\mathcal{C}_1\times...\times \mathcal{C}_d$. For
every subset $\mathcal{E}$ of $\mathbb{R}^d$, and any $\alpha>0$,
introduce the $\alpha$-neighborhood $\mathcal{E}^\alpha$ of
$\mathcal{E}$, namely $\mathcal{E}^\alpha = \{x : \inf_{y\in
\mathcal{E}}|x-y|_{\mathbb{R}^d}\leq \alpha\}$,
$|\cdot|_{\mathbb{R}^d}$ standing for the Euclidean norm on
$\mathbb{R}^d$.\\ We will work under the following regularity
assumptions on $f$ and $f_\ell$, $\ell=1, ..., d$ ($f_\ell$
denoting the density function of $X_\ell$, as in Remark
\ref{rem1}). These functions are supposed to be continuous.
Moreover, we assume the existence of a constant $\alpha>0$ such
that the following assumptions hold.
\begin{tabbing}
$(F.1)\;\;$\=\ $\forall x_\ell\in\mathcal{C}^\alpha_\ell,
f_{\ell}(x_\ell)>0,\
\ell=1, ..., d$,  and  $\forall{\bx}\in\mathcal{C}^\alpha, f({\bx})>0.$  \\
$(F.2)$\>\ $f$ is $k'$-times continuously differentiable on
$\mathcal{C}^\alpha,$ with $k'> kd$.
\end{tabbing}

\noindent Regarding the kernels $K$ and $K_\ell,~ \ell=1, \ldots,
d$, defined in $\mathbb{R}^d$ and $\mathbb{R}$ respectively, they
are assumed to be bounded, integrable to 1, with compact support
and such that,
\begin{tabbing}
$(K.1)\;\;$\=\ $K_\ell$ is of order $k$, $\ell=1, ..., d$.
\end{tabbing}
In addition, we impose the following assumptions on the given
integrating density functions $q_{-\ell}$ and $q_\ell$, $\ell=1,
..., d$.
\begin{tabbing}
$(Q.1)\;\;$ \=\ $q_{-\ell}$  is bounded and
continuous,  $\ell=1, ..., d$.\\
$(Q.2)$\>\ For $\ell=1, ..., d$, $q_\ell$ has a compact support
included in $\mathcal{C}_\ell$ and has \\ \>\ $(k+1)$ continuous
and bounded derivatives  .
\end{tabbing}
Finally, turning our attention to the sequences $(h_n)_{n\geq1}$
and $(h_{\ell,n})_{n\geq1},~\ell=1,...,d$, we will work under the
conditions below.
\begin{tabbing}
$(H.1)\;\;$\=\ $h_n=c'\Big(\frac{\log n}{n}\Big)^{1/(2k'+d)}$,
for a fixed $0<c'<\infty$.\\
$(H.2)$\>\ $h_{\ell,n}= cn^{-1/(2k+1)}$, $\ell=1,...,d,$ for a
fixed $0<c<\infty$.
\end{tabbing}
Some more notation is needed for the statement of our results.
Set, for all ${\bx}\in \mathcal{C}$ and every $\ell=1, ..., d$,
\begin{equation}
b_\ell(x_\ell) = \frac{c^k}{k!}\int_{\mathbb{R}}u^kK_\ell(u)du
\Big((-1)^k
m_\ell^{(k)}(x_\ell)-\int_{\mathbb{R}}m_\ell(z)q_\ell^{(k)}(z)dz
\Big), \label{def_b_l}
\end{equation}
and
\begin{eqnarray}
\sigma_{\ell}^2(x_\ell)=\frac{\int_{\mathbb{R}}K_\ell^2(u)d
u}{cf_\ell(x_\ell)}\int_{\mathbb{R}^{d-1}}H({\bx})\frac{q_{-\ell}^2({\bx}_{-\ell})}{f({\bx}_{-\ell}|x_\ell)}d{\bx}_{-\ell},
\label{def_sigma}
\end{eqnarray}
where
\begin{eqnarray}
H({\bx}) = \mathbb{E}\Big[\frac{\psi^2(Y)}{G(Y)} \big|\ {\bf
X}={\bx}\Big]. \label{def_H}
\end{eqnarray}
We have now all the ingredients to state our main result in
Theorem \ref{th} below.
\begin{theorem}\label{th}
Under the conditions $(C.1$-$2$-$3$-$4)$, $(F.1$-$2)$, $(K.1)$,
$(Q.1$-$2)$ and $(H.1$-$2)$, we have, for every $\ell=1, ..., d$
and all $x_\ell\in\mathcal{C}_\ell$,
\begin{eqnarray}
 \frac{n^{k/(2k+1)}\{\widehat \eta^\star_\ell (x_\ell) - \eta_\ell(x_\ell)\}-
b_\ell(x_\ell)}{\sigma_{\ell}(x_\ell)} \stackrel{\mathcal{L}}{
\longrightarrow} \mathcal{N}(0,1).
\end{eqnarray}
\end{theorem}

\noindent This result naturally implies the following corollary,
which correspond to a refinement of Theorem 3.1 in
\cite{DebbarhViallon3}.

\begin{corollary}\label{cor1}Under the conditions of Theorem \ref{th}, we
have, for every $\ell=1, ..., d$ and all
$x_\ell\in\mathcal{C}_\ell$,
\begin{eqnarray*}
n^{2k/(2k+1)}\E\big(\widehat \eta^\star_\ell(x_\ell) -
\eta_\ell(x_\ell)\big)^2=b_\ell^2(x_\ell) +
\sigma_{\ell}^2(x_\ell)+o(1).
\end{eqnarray*}\end{corollary}

\noindent The proof of Theorem \ref{th} is postponed to Section
\ref{sec_Proofs}. A rough outline of our arguments is as follows.
First, we consider the case where both the density function $f$ of
${\bf X}$ and the function $G$ are known (see Lemma \ref{lem1}
below). In this setting, we show that the estimator defined in
(\ref{estcarbo1}) can be written as an uncensored estimator of the
regression function, and the result of Lemma \ref{lem1} follows
from similar arguments as those developed in \cite{Sarad2000} in
the uncensored setting. Then, from the result of \cite{Ango-Nze},
we can extend Lemma \ref{lem1} to the case where only $G$ is
known. Finally, the uniform consistency of $G^\star_n$ (see, for
instance, \cite{Folder1981}) enables us to conclude the
demonstration of Theorem \ref{th} in the general case.

\section{Application : construction of confidence intervals}\label{sec_Applications}
\setcounter{equation}{0} \setcounter{theo}{0} \setcounter{lem}{0}
\setcounter{cor}{0} \setcounter{rem}{0} \setcounter{prop}{0}
\setcounter{fact}{0}
\subsection{Construction of confidence intervals}
Under the assumption of Theorem \ref{th}, it is straightforward
that, for $\ell=1,...,d$, the interval
\begin{equation*}
[\widehat \eta^\star_\ell (x_\ell) \pm
n^{-k/(2k+1)}(1.96\sigma_{\ell}(x_\ell) - b_\ell(x_\ell))]
\end{equation*}
is a confidence interval for $\eta_\ell (x_\ell)$, at an
asymptotic $95\%$ level, for all $x_\ell\in\mathcal{C}_\ell$.
However, as can be seen in the definition (\ref{def_b_l}), the
bias term $b_\ell$ is a "complex" function involving $m_\ell$ and
$m^{(k)}_\ell$. Bias estimates could be built in by using
estimates of these quantities, but this would result in quite
complex algorithms to derive the confidence bands. Following the
ideas which have been proposed to construct confidence bands for
kernel-type estimators of the regression function (see, e.g.,
Section 4.2 in \cite{Hardle90}), our aim is now to make the bias
term vanish. A close look into the proof of Theorem \ref{th}
reveals that, if the bandwidth $h_{\ell,n}$ is chosen proportional
to $n^{-1/(2k+1)}$ times a sequence that tends slowly to $0$ then
the bias vanishes asymptotically. In this case, for
$\ell=1,...,d$, the interval
\begin{equation*}
[\widehat \eta^\star_\ell (x_\ell) \pm
1.96\sigma_{\ell}(x_\ell)n^{-k/(2k+1)}]
\end{equation*}
provides a confidence interval for $\eta_\ell (x_\ell)$, at an
asymptotic $95\%$ level, for all $x_\ell\in\mathcal{C}_\ell$.
Then, given any consistent estimator $\widehat \sigma_{\ell}$ of
$\sigma_{\ell}$ (making use, for instance, of kernel estimators
for the density functions involved in $\sigma_{\ell}$), we
conclude by Slutsky's Theorem that
\begin{equation}\label{IC}
[\hat \eta^\star_\ell (x_\ell) \pm 1.96\widehat
\sigma_{\ell}(x_\ell)n^{-k/(2k+1)}]
\end{equation}
provides a confidence interval for $\eta_\ell (x_\ell)$, at an
asymptotic $95\%$ level, for all
$x_\ell\in\mathcal{C}_\ell$.\vskip5pt

\subsection{Illustration : a simple simulation study}

In the following paragraph, we present some results from a
simulation study, which especially enables to compare the just
given confidence bands with the almost certainty bands we proposed
in \cite{DebbarhViallon5}.\vskip3pt

\noindent We worked with a sample size $n=1000$, and considered
the case where ${\bf X}=(X_1, X_2)\in\R^2$ (i.e. $d=2$) was such
that $X_1\sim\mathcal{U}(-1,1)$ and $X_2\sim\mathcal{U}(-1,1)$,
where $\mathcal{U}(a,b)$ stands for the uniform law on (a,b). Set
$m_{1}(x)=0.5\times\cos^2(x)$ and $m_2(x)=0.5\times\sin^2(x)$. We
selected $\psi=\1_{\{.\leq 0.9\}}$, and considered the model
$\E[\psi(Y)|X_1=x_1, X_2=x_2]= m_{1}(x_1)+m_{2}(x_2)$. Under this
model, the variable $Y$ was simulated as follows. For each integer
$1\leq i \leq n$, let $p_i=m_{1}(x_{1,i})+ m_{2}(x_{2,i})$ where
$x_{j,i}$ is the $i$-th observed value of the variable $X_j$,
$j=1,2$. Note that $0< p_i< 1$ for every $1\leq i \leq n$. Each
$Y_i$ was then generated as one $\mathcal{U}(0.9-p_i, 1+0.9-p_i)$
variable. Following this proceed ensured that $\mathbb{P}(Y_i\leq
0.9|X_i=x_i)=p_i=m_1(x_{1,i})+ m_2(x_{2,i})$. Regarding the
censoring variable, we generated an i.i.d. sample $C_1, ..., C_n$
such that $C_i \sim\mathcal{U}(0,1)$. This choice yielded, \emph{a
posteriori}, $\p(\delta=1)\simeq 0.2$. We used Epanechnikov
kernels (for $K$, $K_1$ and $K_2$) and selected $q_1=q_2=0.5\times
\1_{[-1,1]}$ (in such a way that the additive component to
estimate were $\eta_{\psi,j}=m_{j} - 0.25$, $j=1,2$). As for the
bandwidth choice, we opted \emph{a priori} for
$h_{1000}=h_{1,1000}=h_{2,1000}=0.2$. \vskip3pt

\begin{figure}[h]
\subfigure[First component]{
\includegraphics[width=0.45\textwidth]{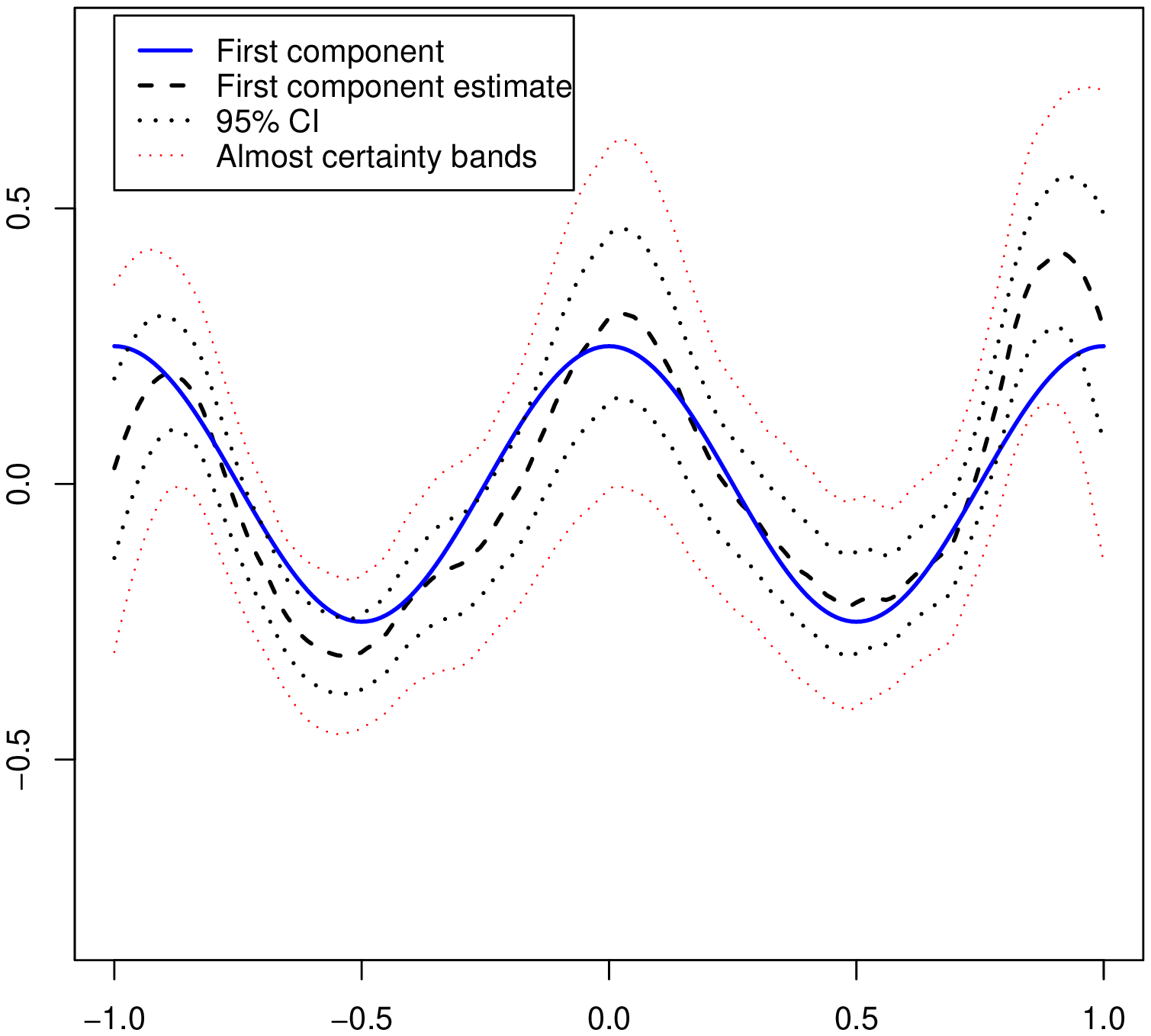}}
\hspace{.3in} \subfigure[Second component]{
\includegraphics[width=0.45\textwidth]{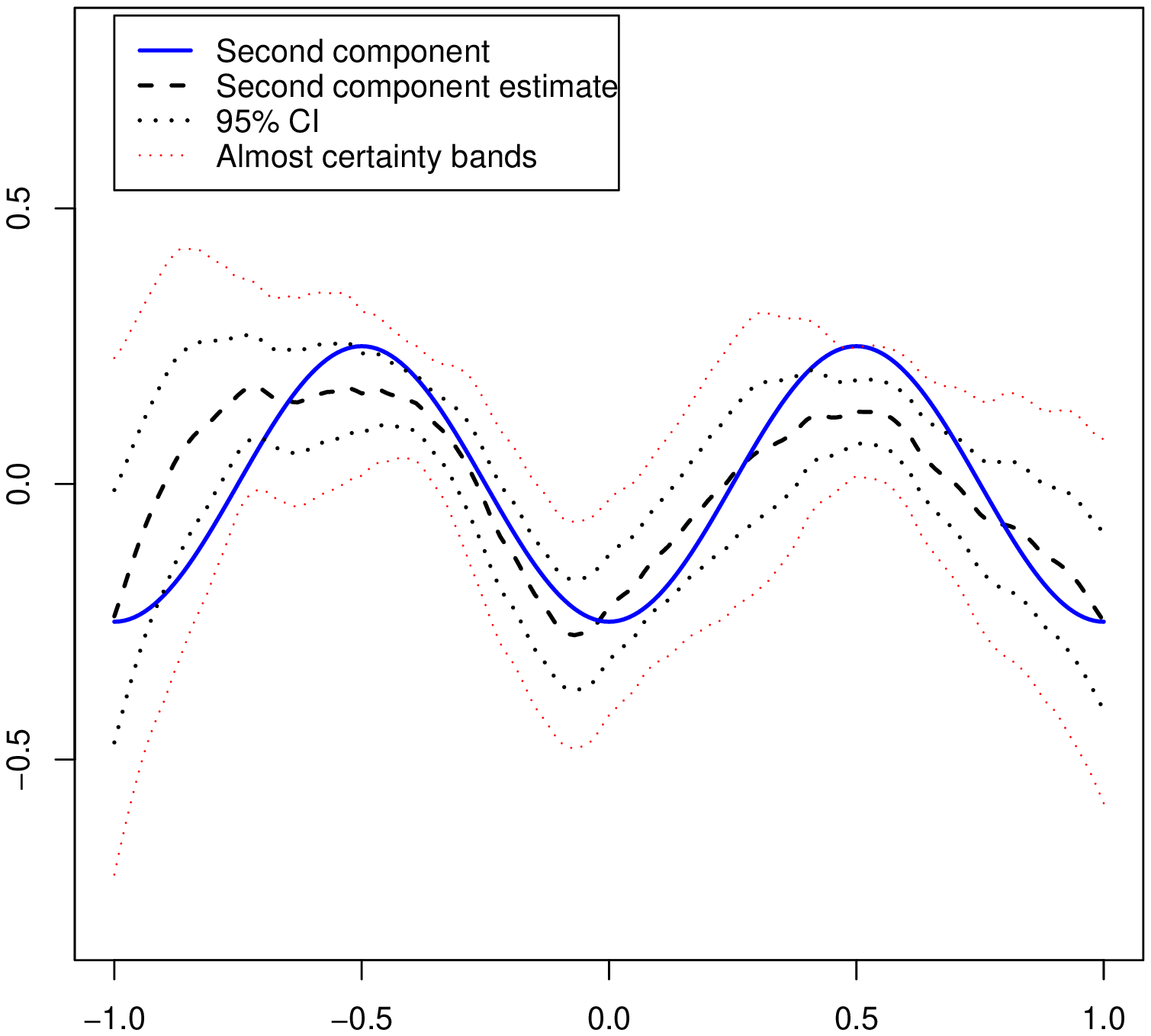}}
\caption{Results of the simulation study for (a) the first
additive component and (b) the second additive component : true
additive components (blue solid line), their estimates (black
dashed line), the 95\% confidence intervals (black dotted lines)
and the almost certainty bands (red dotted lines).}\label{graph}
\end{figure}

\noindent Graphical representations of the results are provided in
Figure \ref{graph}. It can be seen that the confidence intervals
derived from the asymptotic normality are less conservative than
the ones obtained from the uniform law of the logarithm. The price
to pay is however that the true function does not belong to the
former intervals at every $x\in[0,1]$. Therefore, in most
applications, we recommend the construction of both confidence
bands to assess the form of the relationship between the dependant
variable and covariates.

\section{Proof}\label{sec_Proofs}
\setcounter{equation}{0} \setcounter{theo}{0} \setcounter{lem}{0}
\setcounter{cor}{0} \setcounter{rem}{0} \setcounter{prop}{0}
\setcounter{fact}{0}
\subsection{Proof of Theorem \ref{th}}
\subsubsection{The case where both $f$ and $G$ are known}
Recall the definitions (\ref{estcarbo1}) and (\ref{enq1}) of
$\widetilde{m}_{\psi,n}^\star$ and $\widehat \eta^\star_{\ell}$
respectively.  Further denote by $\widetilde{\widetilde
m}_{\psi,n} ({\bx})$ and $\widehat{\widehat
\eta}_{\ell}(x_{\ell})$ the versions of
$\widetilde{m}_{\psi,n}^\star$ and $\widehat \eta^\star_{\ell}$
respectively,  in the case where both $f$ and $G$ are known.
Namely, we have
\begin{equation}
\widetilde {\widetilde m}_{\psi,n} ({\bx}) = \sum_{i=1}^n
\widetilde W_{n,i}({\bx})\frac{\delta_i \psi(Z_i)}{G(Z_i)} ~~
\mbox{with}~ \widetilde W_{n,i}({\bx})= \frac{\prod_{\ell=1}^{d}
\frac{1}{h_{\ell,n}} K_{\ell}\big({\frac{x_{\ell} -
X_{i,l}}{h_{\ell,n}}}\big)}{n f({\bf X}_i)},
\label{estcarbo_fGconnues}
\end{equation}
and
\begin{equation}
\widehat{\widehat \eta}_{\ell}(x_{\ell}) = \int_{\mathbb{R}^{d-1}}
\widetilde {\widetilde m}_{\psi,n}({\bx}) q_{-\ell}({\bx}_{-\ell})
d{\bx}_{-\ell} - \int_{\mathbb{R}^d} \widetilde {\widetilde
m}_{\psi,n}({\bx}) q({\bx}) d{\bx}.\label{def_eta_fGconnues}
\end{equation}
Consider the function $\Psi:\mathbb{R}^2 \longrightarrow
\mathbb{R}$ such that,
\begin{eqnarray}
\Psi(y,c)=\frac{\mathbb{I}_{\{y \leq c\}}\psi(y \wedge
c)}{G(y\wedge c)},~~ \mbox{for all}~ (y,c)\in \mathbb{R}^2.
\label{Psi}
\end{eqnarray}
In view of  (\ref{estcarbo_fGconnues}) and (\ref{Psi}), we have
\begin{eqnarray} \widetilde{\widetilde m}_{\psi,n}({\bx})=
\sum_{i=1}^n \widetilde W_{n,i}({\bx})\Psi(Y_i, C_i).
\label{estcarbo2_fGconnues}
\end{eqnarray}
In the sequel we will make frequent use of a conditional argument,
along with the independence assumption $(C.1)$, that especially
enables us to obtain the following kind of result.
\begin{eqnarray}
m_{\Psi}({\bf X}) := \E(\Psi(Y,C)|{\bf X})  &=& \
\E\Big\{\frac{\mathbb{I}_{\{Y\leq C\}}\psi(Z)} {G(Z)}\Big|\ {\bf
X}\Big\}\nonumber\\
&=&\ \E \Big\{\frac{\psi(Y)}{G(Y)} \E\big[\mathbb{I}_{\{Y\leq
C\}}|{\bf X}, Y\big]\big|\ {\bf X}\Big\} \nonumber\\
&=&\ m_{\psi}({\bf X}).\label{cond_argument}
\end{eqnarray}
Combining this last result with the fact that the quantity
$\Psi(Y_i, C_i)$ is observed (i.e., uncensored) for all
$i=1,...,n$, we have the following particular appealing property :
$\widetilde{\widetilde m}_{\psi,n}({\bx})$ turns out to be an
uncensored estimator of the regression function $m_\psi$. This
property enables to treat the particular case where $G$ is known
with arguments similar to those used in the uncensored
case.\vskip5pt

\noindent We will first establish the following result, which
correspond to Theorem \ref{th} in the case where $f$ and $G$ are
known.
\begin{lemma} \label{lem1} Assume $(C.1$-$2$-$3$-$4)$, $(K.1)$, $(Q.1$-$.2)$
and $(H.2)$ hold. Then, for every $\ell=1,\ldots,d$ and all
$x_\ell\in\mathcal{C}_\ell$,
\begin{eqnarray}
 \frac{n^{k/(2k+1)}\{\widehat {\widehat \eta}_\ell (x_\ell) - \eta_\ell(x_\ell)\}-
b_\ell(x_\ell)}{\sigma_{\ell}(x_\ell)} \stackrel{\mathcal{L}}{
\longrightarrow} \mathcal{N}(0,1).
\end{eqnarray}
\end{lemma}
\noindent {\sc Proof.} In view of the discussion above, we will
mostly borrow the arguments developed in the uncensored case by
\cite{Sarad2000}. Towards this aim, introduce the following
quantities.
\begin{equation}
\tilde \Psi_n( Y_{i},C_i)= \!\Psi (Y_i,C_i)\!\!
\int_{\mathbb{R}^{d-1}}\!\prod_{\ell=2}^{d} \frac{1}{h_{\ell,n}}
K_{\ell}\Big({\frac{x_{\ell} - X_{i,l}}{h_{\ell,n}}}\Big)
\frac{q_{-1}({\bx}_{-1})}{f(X_{i,-1}|X_{i,1})} d{\bx}_{-1},
\label{def_Psi_tilde}
\end{equation}
\begin{equation}
\mathcal{G}({\bf u}_{-1})=\int_{\mathbb{R}^{d-1}}
\prod_{\ell=2}^{d} \frac{1}{h_{\ell,n}}
K_{\ell}\Big({\frac{x_{\ell} - u_{\ell}}{h_{\ell,n}}}\Big)
q_{-1}({\bx}_{-1}) d{\bx}_{-1}, \label{def_G_cal}
\end{equation}
\begin{equation}
\hat \alpha_1(x_1) = \frac{1}{nh_{1,n}} \sum_{i=1}^n \frac{\tilde
\Psi_n(Y_{i},C_i)}{f_1(X_{i,1})} K_{1} \Big(\frac{x_{1}-
X_{i,1}}{h_{1,n}}\Big), \label{def_alpha1}
\end{equation}
\begin{equation}
\tilde m(x_1)= \E\big( \left. \tilde \Psi_n( Y_{i},C_i)\right|
X_{i,1} = x_1\big), \label{def_mtilde}
\end{equation}
\begin{equation}
C_n =  \mu + \int_{\mathbb{R}^{d-1}} \sum_{j=2}^d m_j(u_j)
\mathcal{G}({\bf u}_{-1}) d{\bf u}_{-1},\label{def_Cn}
\end{equation}
\begin{equation}
\hat C_n=\int_{\mathbb{R}^d} \widetilde{\widetilde
m}_{\psi,n}({\bx}) q({\bx})d{\bx}, \label{def_Cn_hat}
\end{equation}
\begin{equation}
C=\int_{\mathbb{R}}m_1(x_1)q_1(x_1)dx_1.
\end{equation}
Next observe that
\begin{eqnarray}
\widehat{\widehat \eta}_1(x_1) - \eta_1(x_1) & = &   \{\hat
\alpha_1(x_1) - \tilde m(x_1)\}+\mathbb{E}\big(\hat C_n - C_n
-C\big), \label{decomp_biais_eta}
\end{eqnarray}
and set
\begin{eqnarray}
&&\beta_1(x_1) = (-1)^k c^k \eta^{(k)}_{1}(x_1) \int_{\mathbb{R}}
\frac{v_1^k}{k!} K_1(v_1)dv_1, \label{beta1}\\
&&\beta_2 = c^k \int_{\mathbb{R}}\Big\{ \int_{\mathbb{R}}
\frac{v_1^k}{k!} m_1^{(k)}(x_1+v_1h_1) K_1(v_1)dv_1 \Big\}q_1(x_1)
dx_1.\label{beta2}
\end{eqnarray}
From $(\ref{decomp_biais_eta})$-$(\ref{beta1})$-$(\ref{beta2})$,
and because of Slutsky's theorem, the proof of Lemma \ref{lem1}
will be completed as soon as the four following results will be
established.
\begin{eqnarray}
&&n^{k/(2k+1)}\ (\hat \alpha_1(x_1)-\tilde m(x_1) -\beta_1(x_1))
\stackrel{\mathcal{L}}{ \longrightarrow}
\mathcal{N}(0,\sigma_{1}^{2}(x_{1})), \label{decomp_lem1_1}\\
&& n^{k/(2k+1)}\ \mathbb{E}\big(\hat C_n-C_n-C\big)= \beta_2 + o(1), \label{decomp_lem1_2}\\
&&n^{2k/(2k+1)}\ {\rm Var} (\hat C_n)=o(1)\label{decomp_lem1_3},\\
&& n^{2k/(2k+1)}\ {\rm Cov}(\hat C_n, \hat\alpha_1(x_1))=o(1).
\label{decomp_lem1_4}
\end{eqnarray}
{\em Proof of (\ref{decomp_lem1_1}):} In a first step, our aim is
to show that
\begin{eqnarray}
nh_{1,n} \{\hat \alpha_1 (x_1) - \E(\hat \alpha_1 (x_1)\}
\longrightarrow \mathcal{N}\big(0, \sigma_1^2(x_1)\big).
\label{norm1}
\end{eqnarray}
We claim that
\begin{eqnarray}
nh_{1,n} {\rm Var}( \hat \alpha_{1}(x_{1}))\rightarrow
\sigma_{1}^{2}(x_{1}) \mbox{ as } n\rightarrow\infty,
\label{var_alpha}
\end{eqnarray}
where $\sigma_{1}^{2}(x_{1})$ is as in (\ref{def_sigma}).
Recalling (\ref{def_alpha1}), note that
\begin{eqnarray*}
{\rm Var}( \hat \alpha_{1}(x_{1}))& = & \frac{1}{n^2}
\sum_{i=1}^{n} \mathbb{E}\Big\{ \frac{\tilde \Psi_n(
Y_{i},C_i)}{f_{1}(X_{i,1})h_{1,n}} K_{1}\Big(\frac {x_{1} -
X_{i,1}}{h_{1,n}}\Big)\Big\}^2 \\
&&  - \frac{1}{n^2} \sum_{i=1}^{n} \Big\{
\mathbb{E}\Big(\frac{\tilde \Psi_n(
Y_{i},C_i)}{f_{1}(X_{i,1})h_{1,n}} K_{1}\Big(\frac {x_{1}
- X_{i,1}}{h_{1,n}}\Big)\Big)\Big\}^2 \\
& = & \frac{1}{n} \left[ \frac{1}{h_{1,n}} \Phi_{1,n}(x_{1}) -  [
\Gamma_{1,n}(x_{1})]^2 \right],
\end{eqnarray*}
with
\begin{eqnarray*}
\Gamma_{1,n}(x_{1}) & = & \mathbb{E} \Big\{ \frac{1}{h_{1,n}}
K_{1}\Big(\frac{x_{1} - X_{i,1}}{h_{1,n}}\Big) \frac{\tilde
\Psi_n( Y_{i},C_i)}{f_{1}(X_{i,1})} \Big\}, \\
\Phi_{1,n}(x_{1})  & = & \mathbb{E}\Big\{\frac{1}{h_{1,n}}
K_{1}^{2}\Big(\frac{x_{1} - X_{i,1}}{h_{1,n}}\Big) \frac {\tilde
\Psi_n^2( Y_{i},C_i)}{f_{1}^{2}(X_{i,1})} \Big\}.
\end{eqnarray*}
Using classical conditioning arguments and recalling the
definition (\ref{def_mtilde}) of $\tilde m$, it is straightforward
that
\begin{eqnarray*}
&&\Gamma_{1,n}(x_{1}) =  \int_{\mathbb{R}} \frac{1}{h_{1,n}}
K_{1}\left(\frac{x_{1} - u_{1}}{h_{1,n}}\right) \tilde m(u_{1})
du_{1}, \\
&&\Phi_{1,n}(x_{1}) =  \int_{\mathbb{R}} \frac{1}{h_{1,n}}
K_{1}^{2}\left(\frac{x_{1} - u_{1}}{h_{1,n}}\right) \frac
{\mathbb{E}(\tilde \Psi_n^2( Y_{i},C_i) \mid X_{i,1}=
u_{1})}{f_{1} (u_{1})}\ du_{1}.
\end{eqnarray*}
In other respect, the definitions (\ref{def_Psi_tilde}) and
(\ref{def_mtilde}), when combined with the argument used in
(\ref{cond_argument}), yield
\begin{eqnarray*}
| \tilde m(u_{1}) | & = &  | \mathbb{E}( \tilde \Psi_n(Y_{i},C_i) \mid X_{i,1} = u_{1}) | \\
& = &  \Big|\mathbb{E}\Big( \frac{\delta \psi(Z_i)}{G(Z_i)}
\frac{{\mathcal G}({\bf X}_{i,-1})}{f({\bf X}_{i,-1}|\
X_{i,1})}\big| X_{i,1} = u_{1} \Big) \Big|\\
& = & \big| \int_{\mathbb{R}^{d-1}} m_\psi({\bf u} ) {\mathcal
G}({\bf u}_{-1})  d{\bf u}_{-1} \big| <\infty,
\end{eqnarray*}
in such a way that
\begin{eqnarray}
\frac{1}{n}\ \Gamma_{1,n}^{2}(x_{1}) \rightarrow 0\mbox{ as }
n\rightarrow \infty. \label{controle_gamma1}
\end{eqnarray}
Moreover, using once again the  argument used to derive
(\ref{cond_argument}) and keeping in mind the definition
(\ref{def_H}) of $H$, it is easy to derive that, under $(C.1)$,
\begin{eqnarray}
&&\mathbb{E}(\tilde \Psi_n^2( Y_{i},C_i) \mid X_{i,1} = x_{1})\nonumber \\
& =&\mathbb{E}\Big\{ \Big( \frac{\delta
\psi(Z_i)}{G(Z_i)}\Big)^{2} \Big(\frac{{\mathcal G}({\bf
X}_{i,-1})}{f({\bf X}_{i,-1} \mid
X_{i,1})}\Big)^2 \Big|\ X_{i,1}=x_{1} \Big\}, \nonumber \\
&=&\mathbb{E}\Big\{\mathbb{E}\Big[\Big(\frac{\delta
\psi(Z_i)}{G(Z_i)}\Big)^{2} \Big| X_{i}\Big] \frac{ {\mathcal
G}^2({\bf X}_{i,-1})}{ f^{2}({\bf X}_{i,-1}|
X_{i,1})} \Big|\ X_{i,1} = x_{1} \Big\},\nonumber \\
&= & \int_{\mathbb{R}^{d-1}} H({\bf u}) \frac{{\mathcal G}^2({\bf
u}_{-1})}{f^2({\bf u}_{-1} | x_{1})}\
f({\bf u}_{-1}\mid x_{1}) d{\bf u}_{-1}, \nonumber \\
& = & \int_{\mathbb{R}^{d-1}} H({\bf u}) \frac{{\mathcal G}^2({\bf
u}_{-1})}{f({\bf u}_{-1} \mid x_{1})}\ d{\bf
u}_{-1}.\label{decomp_var1}
\end{eqnarray}
Next, making use of the classical change of variable
$v_{\ell}h_{\ell,n} = u_{\ell} - x_{\ell}$ along with a Taylor
expansion of order $k$ (which is rendered possible by $(Q.2)$), we
readily have by $(K.1)$, for a given $0 < \theta < 1$,
\begin{eqnarray}
&     & \int_{\mathbb{R}^{d-1}} \prod_{\ell=2}^{d}
\frac{1}{h_{\ell,n}} K_{\ell}\left(\frac{x_{\ell} -
u_{\ell}}{h_{\ell,n}}\right) q_{\ell}( x_{\ell})
d{\bx}_{-1} - q_{-1}({\bf u}_{-1})\nonumber\\
& =   & \int_{\mathbb{R}^{d-1}} \prod_{\ell=2}^{d}  \left(
K_{\ell}(v_{\ell}) \left[\frac{v_{\ell}^{k}h_{\ell,n}^{k}}{k!}\
q_{\ell}^{(k)}
(\theta v_{\ell} h_{\ell,n} + u_{\ell}) \right]\right) d{\bf v}_{-1} \nonumber\\
&  =  & o(1)\label{decomp_var2},
\end{eqnarray}
Combining (\ref{decomp_var1}) and (\ref{decomp_var2}), we get
\begin{equation}
\mathbb{E}( \tilde \Psi_n^2(Y_i,C_i) \mid X_{i,1} = x_{1} ) =
\int_{\mathbb{R}^{d-1}} H({\bx}) \frac{q_{-1}^{2}
({\bx}_{-1})}{f({\bx}_{-1} \mid x_{1} )} ~d{\bx}_{-1} + o(1).
\label{EPsi_n}
\end{equation}
In addition, setting $ \Phi (x_{1}) = \displaystyle
\int_{\mathbb{R}^{d-1}} H({\bx}) \frac{q_{-1}^{2}
({\bx}_{-1})}{f({\bx}_{-1} \mid x_{1} )} ~d{\bx}_{-1}$ and using
once again the change of variable $ v_{1} h_{1,n} = x_{1}- u_{1}
$, we obtain
\begin{eqnarray*}
\Phi_{1,n} (x_{1}) & = & \int_{\mathbb{R}} \frac{K_{1}^{2}
(v_{1})}{ f_{1}(x_{1} -h_{1} v_{1})} \mathbb{E}(\tilde
\Psi^2_n(Y_{i},C_i)|\ X_{i,1} = x_{1} - h_{1} v_{1}) dv_{1} \\
& = & \int_{\mathbb{R}} K_{1}^{2} (v_{1}) \bigg( \frac{
\mathbb{E}(\tilde \Psi^2_n(Y_{i},C_i) \mid X_{i,1} = x_{1} - h_{1}
v_{1})}{ f_{1}( x_{1}-h_{1} v_{1})}
-\frac{\Phi(x_{1})}{f_{1}(x_{1})} \bigg)
dv_{1}\\
&   &  + \frac{\Phi(x_{1})}{f_{1}(x_{1})} \int_{\mathbb{R}}
K_1^{2} (v_{1} ) dv_{1}.
\end{eqnarray*}
But, by $(C.3),(K.1), (F.1)$ and $(Q.1$-$2)$, it is easily shown
that the quantity $|\E(\tilde \Psi(Y_i,C_i)|\ X_{i,1} = u_1)/
f_{1}(u_{1}) -\Phi(x_{1})/f_{1}(x_{1})|$ is bounded. Therefore,
(\ref{EPsi_n}) when combined with Lebesgue's dominated convergence
Theorem enables us to conclude that
\begin{eqnarray}
\Phi_{1,n} (x_{1}) \rightarrow
\frac{\Phi(x_{1})}{f_{1}(x_{1})}\int_{\mathbb{R}} K_1^{2}(v_{1})
dv_{1}. \label{controle_phi}
\end{eqnarray}
From (\ref{controle_gamma1}) and (\ref{controle_phi}), the claim
(\ref{var_alpha}) is proved. \vskip5pt

\noindent Now, we set
\begin{eqnarray}
& &\widetilde T_{i,n} =  \frac{1}{h_{1,n}} \frac{\tilde
\Psi_n(Y_{i},C_i)}{f_1(X_{i,1})} K_{1} \Big(\frac{x_{1}-
X_{i,1}}{h_{1,n}}\Big), \quad
T_{i,n} = \widetilde T_{i,n} - \mathbb{E} \widetilde T_{i,n}, \nonumber\\
&&\mbox{and} \quad s_n^2 = \mbox{Var}(T_{i,n})= n \mbox{Var}(\hat
\alpha_1(x_1)). \label{def_sn}
\end{eqnarray}
For all $\varepsilon > 0$, we have
\begin{eqnarray}
\mathbb{E}\Big\{\frac{T^2_{i,n}}{ns^2_n}
\mathbb{I}_{\{|\frac{T_i}{\sqrt{n}S_n}|\geq~ \varepsilon \}}\Big\}
& \leq & \frac{M_1}{h_{1,n}^2n
s_n^2}\mathbb{P}\Big(\big|\frac{T_{i,n}}{\sqrt{n}S_n}\big|\geq \varepsilon\Big)\nonumber\\
& \leq &\frac{M_1}{h_{1,n}^2n s_n^2}
\frac{\mathbb{E}(T_{i,n}^2)}{\varepsilon^2n s_n^2}\nonumber\\
& \leq &\frac{M_1\mathbb{E}(T_{i,n}^2)}{(\sigma_1^2(x_1)+o(1))^2\varepsilon^2n^4}\nonumber\\
& \leq & \frac{M_2}{\varepsilon^2 h_{1,n}^2 n^4},
\label{cond_normality}
\end{eqnarray}
where $M_1$ and $M_2$ are two finite and positive constants.
Combining (\ref{cond_normality}) with the fact that
$T_{i,n}/s_n^2\rightarrow0$ (which follows from (\ref{var_alpha})
and (\ref{def_sn})), we can apply the normal convergence criterion
(see, e.g. \cite{Loeve}, p.295) to obtain
\begin{eqnarray}
\frac{1}{\sqrt{n}s_n}\sum_{i=1}^n
\big(T_{i,n}-\mathbb{E}T_{i,n}\big)\stackrel{\mathcal{L}}{
\longrightarrow} \mathcal{N}(0,1). \label{normality1}
\end{eqnarray}
Finally, (\ref{norm1}) readily comes from (\ref{def_alpha1}),
(\ref{var_alpha}), (\ref{def_sn}) and (\ref{normality1}).
\vskip5pt

\noindent Now, our aim is to evaluate the term $|\tilde m(x_1) -
\mathbb{E}\hat \alpha_1(x_1)|$. First, from (\ref{def_Psi_tilde}),
(\ref{def_G_cal}) and (\ref{def_mtilde}), note that
\begin{eqnarray*}
\tilde m(x_1)&=&\mathbb{E}\Big( \tilde \Psi_n( Y_{i},C_i)\big|\
X_{i,1} =
x_1 \Big) \\
&=& \mathbb{E} \Big( \frac{\Psi(Y_i, C_i)}{{f({\bf
X}_{i,-1}|X_{i,1})}}\ \mathcal{G}({\bf X}_{i,-1}) \big|\
X_{i,1}=x_1 \Big).
\end{eqnarray*}
Then, using a conditioning argument along with the independence
assumption $(C.1)$, we get
\begin{eqnarray*}
\tilde m(x_1)&=&\mathbb{E} \Big\{ \frac{\mathbb{E}(\Psi(Y_i,
C_i)|{\bf X}_i)}{{f({\bf X}_{i,-1}|X_{i,1})}}\ \mathcal{G}({\bf
X}_{i,-1})
\big|\ X_{i,1}=x_1 \Big\} \\
&=& \mathbb{E} \Big\{ \frac{\mathbb{E}(\psi(Y_i)|{\bf
X}_i)}{{f({\bf X}_{i,-1}|X_{i,1})}}\ \mathcal{G}({\bf X}_{i,-1})
\big|\ X_{i,1}=x_1 \Big\} \\
&=& \int_{\mathbb{R}^{d-1}} m_\psi(x_1,{\bf
u}_{-1})\mathcal{G}({\bf u}_{-1})\ d{\bf u}_{-1}.
\end{eqnarray*}
Thus, by $(K.1)$ and $(C.4)$, a Taylor expansion yields
\begin{eqnarray*}
&& \mathbb{E}(\hat \alpha_1(x_1) ) - \tilde m(x_1)\\
&=& \int_{\mathbb{R}}\frac{1}{h_{1,n}}
\tilde m(u_1) K_1\Big( \frac{x_1-u_1}{h_{1,n}} \Big) du_1 -\tilde m(x_1)\\
&=& \int_{\mathbb{R}} \big[\tilde m(x_1-v_1h_{1,n})-\tilde
m(x_1)\big] K_1(v_1)dv_1\\
&=&\int_{\mathbb{R}} \int_{\mathbb{R}^{d-1}}
\big[m_\psi(x_1-v_1h_{1,n},{\bf u}_{-1}) - m_\psi(x_1,{\bf
u}_{-1})\big] \mathcal{G}({\bf u}_{-1}) d{\bf u}_{-1} K_1(v_1)dv_1\\
& = & \int_{\mathbb{R}} \int_{\mathbb{R}^{d-1}}
\left[\frac{(-h_{1,n}v_1)^k}{k!} \frac{\partial^k
m_{\psi}}{\partial x_1^k}(x_1-h_{1,n} v_1,{\bf u}_{-1} )\right]
\mathcal{G}({\bf u}_{-1}) d{\bf u}_{-1}  K_1(v_1)dv_1\\
&& + o(h_{1,n}^k)\\
& = &  h_{1,n}^k\beta_1(x_1) + o(h_{1,n}^k).
\end{eqnarray*}
By combining this last result with (\ref{norm1}), we conclude to
(\ref{decomp_lem1_1}).$\sqcup\!\!\!\!\sqcap$ \vskip10pt

\noindent \emph{Proof of (\ref{decomp_lem1_2})}. Keep in mind the
definitions (\ref{def_Cn}) and (\ref{def_Cn_hat}) of $C_n$ and
$\hat C_n$. Then, according to Fubini's Theorem and under the
additive model assumption,
\begin{eqnarray*}
\mathbb{E}(\hat C_n -C_n) & = &
\mathbb{E}\Big\{\int_{\mathbb{R}^d} \widetilde {\widetilde
m}_{\psi,n}({\bx}) q({\bx})d{\bx} - \mu - \int_{\mathbb{R}^{d-1}}
\sum_{j=2}^d m_j(u_j)\mathcal{G}({\bf
u}_{-1})d{\bf u}_{-1}\Big\}\\
&& - \int_{\mathbb{R}^{d-1}} \sum_{j=2}^d m_j(u_j)
\mathcal{G}({\bf u}_{-1}) d{\bf u}_{-1}\\
& = & \sum_{j=1}^d \int_{\mathbb{R}^d}
\frac{1}{h_{1,n}}m_{j}(u_j)\mathcal{G}({\bf
u}_{-1})\int_{\mathbb{R}}K_1\Big(\frac{x_1-
u_1}{h_{1,n}}\Big)q_1(x_1)dx_1 d{\bf u}\\
&&-\int_{\mathbb{R}^{d-1}} \sum_{j=2}^d m_j(u_j) \mathcal{G}({\bf
u}_{-1}) d{\bf u}_{-1}\\
&  = &\int_{\mathbb{R}}\int_{\mathbb{R}}
\frac{1}{h_{1,n}}m_1(u_1)K_1\Big(\frac{x_1-
u_1}{h_{1,n}}\Big)q_1(x_1) dx_1du_1.
\end{eqnarray*}
But, by $(C.4)$ and $(K.1)$ and using a Taylor expansion, we get,
\begin{eqnarray}
&   & \mathbb{E}(\hat C_n -C_n) - C \nonumber\\
& = &
\int_{\mathbb{R}}\int_{\mathbb{R}}q_1(x_1)m_1(x_1+h_{1,n}v_1)K_1(v_1)dv_1dx_1 - C\nonumber\\
& = & \int_{\mathbb{R}}\int_{\mathbb{R}}q_1(x_1)
[m_1(x_1+h_{1,n}v_1)-m_1(x_1)]K_1(v_1)dv_1dx_1\nonumber\\
& = &
\int_{\mathbb{R}}\int_{\mathbb{R}}q_1(x_1)\Big[\frac{h_{1,n}^k
v_1^k}{k!} m_1^{(k)}(x_1+v_1h_1)\Big]K_1(v_1)dv_1dx_1 + o(h_{1,n}^k)\nonumber\\
& =& h_{1,n}^k \beta_2 + o(h_{1,n}^k), \label{res_decomp2_2}
\end{eqnarray}
which allows us conclude to
(\ref{decomp_lem1_2}).$\sqcup\!\!\!\!\sqcap$ \vskip10pt

\noindent \emph{Proof of (\ref{decomp_lem1_3})}. In view of the
definitions (\ref{estcarbo_fGconnues}), (\ref{def_Psi_tilde}),
(\ref{def_G_cal}) and  (\ref{def_Cn_hat}), using the boundedness
of ${\mathcal G}$ (which is ensured by $(Q.1$-$2)$) and $\Psi$
(which is ensured by $(C.3)$) along with the fact that, by
$(F.1)$, $f$ is bounded away from 0, we have, for a given $M_3>0$,
\begin{equation}
{\rm Var} (\hat C_n)
\leq
\frac{M_3}{nh_{1,n}^2}\mathbb{E}\Big(\int_{\mathbb{R}}K_1\Big(\frac{x_1-X_{i,1}}{h_{1,n}}\Big)q_1(x_1)dx_1\Big)^2
= \mathcal{O}\Big(\frac{1}{n}\Big), \label{varCn}
\end{equation}
which naturally implies
(\ref{decomp_lem1_3}).$\sqcup\!\!\!\!\sqcap$ \vskip10pt

\noindent \emph{Proof of (\ref{decomp_lem1_4})} Using the
Cauchy-Schwartz inequality, we infer from (\ref{var_alpha}) and
(\ref{varCn}) that, \[{\rm Cov}(\hat C_n, \hat\alpha_1(x_1))\leq
({\rm Var}\ \hat C_n)^{1/2}({\rm Var}\ \hat\alpha_1(x_1))^{1/2} =
\mathcal{O}\big(n^{-1/(2k+1)}\big),\] which implies
(\ref{decomp_lem1_4}).$\sqcup\!\!\!\!\sqcap$ \vskip8pt

\noindent As already pointed out, the proof of Lemma $\ref{lem1}$
is readily completed by combining (\ref{decomp_lem1_1}),
(\ref{decomp_lem1_2}), (\ref{decomp_lem1_3}) and
(\ref{decomp_lem1_4}). $\sqcup\!\!\!\!\sqcap$

\subsubsection{The case where $f$ is unknown but $G$ is known}
The key idea in this case is to use the uniform consistency of
$\hat f_n$ (see, e.g., \cite{Ango-Nze}) along with the following
decomposition,
\begin{equation}\label{decomposition}
\frac{1}{\hat f_n}  =  \frac{1}{f} - \frac{\hat f_n - f}{\hat f_n
f}.
\end{equation}
When  the density $f$ is unknown and $G$ is known, the additive
components estimates are defined as follows, for
$\ell=1,\ldots,d,$
\begin{eqnarray}
{\widehat \eta}_{\ell}(x_{\ell}) = \int_{\mathbb{R}^{d-1}}
{\widetilde{m}}_{\psi,n}({\bx}) q_{-\ell}({\bx}_{-\ell})
d{\bx}_{-\ell} - \int_{\mathbb{R}^d}{\widetilde m}_{\psi,n}
({\bx}) q({\bx}) d{\bx}. \label{eta_l_G_connu}
\end{eqnarray}
where
\begin{eqnarray}
{\widetilde m}_{\psi,n} ({\bx}) = \sum_{i=1}^n
W_{n,i}({\bx})\frac{\delta_i \psi(Z_i)}{G(Z_i)}.
\label{est_Gconnu}
\end{eqnarray}
We will establish the following result.
\begin{lemma}\label{lemmeapprox1} Under the hypotheses of Theorem
\ref{th}, we have
\begin{eqnarray}\label{Gc}
\sup_{x_\ell \in {\mathcal C}_\ell} |{\widehat
\eta}_{\ell}(x_{\ell}) - \widehat {\widehat
\eta}_{\ell}(x_{\ell})| =\mathcal{O}\left(\sqrt{\frac{\log
n}{nh_n^d}}\right)~~\mbox{a.s. .}
\end{eqnarray}
\end{lemma}\vskip5pt \noindent {\sc Proof.} First note that the term
$(nh_{1,n}^d)^{-1}\sum_{i=1}^n |K(({\bx}-{\bf X}_i)/h_{1,n})|$ is
almost surely uniformly bounded on $\mathcal{C}$ under the
assumptions we made on $K$ and $f$. Moreover, $f$ and then $f_n$
(for $n$ large enough) are bounded away from 0 (see $(F.1)$).
Then, in view of the definitions (\ref{estcarbo1}) and
(\ref{estcarbo_fGconnues}), along with the decomposition
(\ref{decomposition}), we get, by $(F.2)$ and $(H.1)$, that, for a
given $C_1>0$,
\begin{eqnarray}
\sup_{{\bx}\in \mathcal{C}}\sum_{i=1}^n|W_{i,n}({\bx})-
\widetilde W_{i,n}({\bx})| &\leq& C_1 \sup_{{\bx}\in \mathcal{C}}|\hat f_n({\bx})-f({\bx})| \nonumber\\
& = & \mathcal{O}\left(\sqrt{\frac{\log
n}{nh_n^d}}\right)~~\mbox{a.s.},
\end{eqnarray}
where we used the following result, due to \cite{Ango-Nze},
\begin{eqnarray*}
\sup_{{\bx} \in \mathcal{C}}|\hat f_n({\bx}) - f({\bx})| =
\mathcal{O}\left(\sqrt{\frac{\log n}{nh_n^d}}\right)~~\mbox{a.s. ,
under $(H.1)$ and $(F.2)$.}
\end{eqnarray*}
Next, under the assumptions $({\bf A}), (C.2)$ and $(C.3)$, we
have $\max_{i} \psi(Z_i)/G(Z_i)$ $<\infty$. Thus, from
(\ref{def_eta_fGconnues}), (\ref{eta_l_G_connu}) and
(\ref{est_Gconnu}), we conclude that, for a given $C_2>0$,
\begin{eqnarray*}
\sup_{x_\ell \in {\mathcal C}_\ell} |{\widehat
\eta}_{\ell}(x_{\ell}) - \widehat {\widehat
\eta}_{\ell}(x_{\ell})| &\leq & 2 \sup_{{\bx}\in \mathcal{C}}
\big|\widetilde {\widetilde
m}_{\psi,n}({\bx})- {\widetilde m}_{\psi,n}({\bx})\big| \\
&\leq & 2C_2 \sup_{{\bx}\in \mathcal{C}}\sum_{i=1}^n|W_{i,n}({\bx})- \widetilde W_{i,n}({\bx})|\\
& =&\mathcal{O}\left(\sqrt{\frac{\log
n}{nh_n^d}}\right)~~\mbox{a.s.,}
\end{eqnarray*}
which is Lemma \ref{lemmeapprox1}.$\sqcup\!\!\!\!\sqcap$
\vskip10pt

\subsubsection{The case where both $f$
and $G$ are unknown} We have the following decomposition.
\begin{eqnarray*}
\sup_{x_\ell\in \mathcal{C}_\ell}| \widehat \eta_\ell^\star
(x_\ell)- \eta_\ell(x_\ell)| &\leq & \sup_{x_\ell\in
\mathcal{C}_\ell}|\widehat \eta_\ell^\star (x_\ell)-\widehat
\eta_\ell (x_\ell)| + \sup_{x_\ell\in
\mathcal{C}_\ell}|\widehat \eta_\ell (x_\ell)- \eta_\ell (x_\ell)|
 \end{eqnarray*}

\begin{lemma}\label{lemmeapprox2} Under the assumptions of Theorem
\ref{th}, we have
\begin{equation}
\sup_{x_\ell\in \mathcal{C}_\ell}|\widehat \eta_\ell^\star
(x_\ell)-\widehat \eta_\ell (x_\ell)| = o\big(n^{-k/(2k+1)}\big)
\mbox{ a.s. }.
\end{equation}
\end{lemma} \vskip5pt \noindent {\sc Proof.} Observe that
\begin{eqnarray}
\sup_{x_\ell\in
\mathcal{C}_\ell}|\widehat{\eta}^\star_{\ell}(x_\ell)
-\widehat{\eta}_\ell(x_\ell)| \leq 2 \sup_{{\bx}\in
\mathcal{C}}|\widetilde{m}_{\psi,n}^\star({\bx})-\widetilde{m}_{\psi,n}({\bx})|.
\label{fin1}
\end{eqnarray}
First consider the case where $({\bf A})(i)$ holds. Under the
assumptions of Theorem \ref{th}, we have
\begin{eqnarray}
&&
|\widetilde{m}_{\psi,n}^\star({\bx})-\widetilde{m}_{\psi,n}({\bx})|\nonumber\\
&\leq & M \sum_{i=1}^n |W_{n,i}({\bx})| \sup_{y\leq \tau_0}
\big|\frac{1}{G(y)}-\frac{1}{G^\star_n(y)}\big| \nonumber \\
&\leq & M \sum_{i=1}^n |W_{n,i}({\bx})| \sup_{y\leq \tau_0}
|G(y)-G^\star_n(y)|\sup_{y \leq \tau_0}\frac{1}{G^\star_n(y)G(y)}
, \label{pb1}
\end{eqnarray}
where $M$ is as in $(C.3)$. Since $\tau_0<T_H$, the iterated law
of the logarithm of \cite{Folder1981} ensures that
\begin{equation}\label{foldes2}
\sup_{y \leq \tau}
|G(y)-G^\star_n(y)|=\mathcal{O}\Big(\sqrt{\frac{\log \log
n}{n}}\Big) \mbox{ a.s. }.
\end{equation}
Besides, by the conditions imposed on $K$ and $f$, the term
$\sum_{i=1}^n |W_{n,i}({\bx})|$ is almost surely uniformly
bounded. Combining this last result with (\ref{foldes2}), it
follows that
\begin{eqnarray}
\sup_{{\bx}\in\mathcal{C}}|\widetilde{m}_{\psi,n}^\star({\bx})-\widetilde{m}_{\psi,n}({\bx})|=\mathcal{O}\Big(\sqrt{\frac{\log
\log n}{n}}\Big) \mbox{ a.s. .} \label{fin2}
\end{eqnarray}
From (\ref{fin1}) and (\ref{fin2}), we readily conclude to the
result of Lemma \ref{lemmeapprox2} in the case where $({\bf
A})(i)$ holds. \vskip3pt

\noindent In the case where $({\bf A})(ii)$ holds, the proof
follows from the same lines as above, making use of either the
iterated law of the logarithm of \cite{GuLai} (if $({\bf A})(ii)$
holds with $p=1/2$) or Theorem 2.1 of \cite{ChenLo} (if $({\bf
A})(ii)$ holds with $k/(2k+1)<p<1/2$) instead of the iterated law
of the logarithm of \cite{foldesrejto}. The details are omitted.
$\sqcup\!\!\!\!\sqcap$\vskip5pt

\vskip10pt Finally, putting the results of Lemma \ref{lem1}, Lemma
\ref{lemmeapprox1} and Lemma \ref{lemmeapprox2} all together
achieves the demonstration of Theorem \ref{th}.


\end{document}